\documentclass{amsart}
\usepackage{amsmath}
\usepackage{amsfonts}
\usepackage{amscd}
\usepackage{enumerate}
\usepackage{fancyheadings}
\usepackage{bbold}
\usepackage{amsthm}
\newtheorem{thm}{Theorem}[section]
\newtheorem{cor}[thm]{Corollary}
\newtheorem{conj}[thm]{Conjecture}

\newtheorem{lem}[thm]{Lemma}
\newtheorem{fact}[thm]{Fact}

\theoremstyle{definition}

\begin{document}

\title[Eigenvalues of Words in Two Positive Definite Letters]{Eigenvalues of Words in Two Positive Definite Letters}

\author{Christopher J. Hillar}%
\address{Mathematics Department, University of California, Berkeley, CA 94720}%
\email{chillar@math.berkeley. edu}%

\author{Charles R. Johnson}%
\address{Mathematics Department, College of William and Mary, Williamsburg, VA
23187-8795}%
\email{crjohnso@math.wm.edu}%

\thanks{This research was conducted, in part, during the summer of 1999 at
the College of William and Mary's Research Experiences for
Undergraduates program}

\begin{abstract}
The question of whether all words in two real positive definite letters
have only positive eigenvalues is addressed and settled (negatively).
This question was raised some time ago in connection with a
long-standing problem in theoretical physics. A large class of words that do
guarantee positive eigenvalues is identified, and considerable evidence
is given for the conjecture that no other words do. In the process, a
fundamental question about solvability of symmetric word equations is
encountered.
\end{abstract}

\maketitle

\section{Introduction}
A \textit{word} is a juxtaposed sequence of letters chosen (with
repetition allowed) from a given alphabet. We shall be concerned here
with an alphabet of two letters, $\{A, B\}$, so that a sample word would
be \textit{AABABBBAAB}; thus, hereafter ``word'' means one over a two-letter
alphabet. The \textit{length} of a word is the total number of letters
present (including repetitions); the sample word has length 10. We shall
be interested in the combinatorial structure of words as abstract
objects, but, often, we will interpret a word as the matrix resulting
from the substitution of two independent positive definite matrices for
$A$ and $B$. The eigenvalues and trace of the resulting matrix will be
our primary interest.

The initial motivation comes from a chain of three questions raised by
Lieb \cite{L}, stemming from issues in quantum physics \cite{BMV}. In
addition Pierce raised Question~3 below from an independent source
\cite{P}. The three questions are the following:

\textit{Question} 1. Does the polynomial $p(t)$, defined by
$p(t)= \mathrm{Tr}[(A+Bt)^m]$, have all positive coefficients whenever
$A$ and $B$ are positive definite matrices?

Since the coefficient of $t^k$ in $p(t)$ is the trace of the sum of all
words in $A$ and $B$ with length $m$ and $k\,B$'s, the following, which
could help answer Question 1, has also been asked \cite{L}.

\textit{Question} 2. Is the trace of a given word positive for all
positive definite $A$ and $B$?

Since a matrix with positive eigenvalues necessarily has positive trace, a
yet more precise question has also been raised \cite{L}, \cite{P}.

\textit{Question} 3. Are all the eigenvalues of a given word positive
for all positive definite $A$ and $B$?

In addition, these particular questions and a number of natural issues
they raise seem central to matrix analysis. Since we became interested
in them (thanks to Lieb and Pierce), we have learned that a number
of different investigators (including us) have tested them empirically
by trying many different words and calculating the eigenvalues for many
(tens of thousands) different randomly generated pairs of matrices of
different sizes. To our knowledge, no one turned up a counterexample via
such simulation, rendering Question 3 all the more interesting. Indeed,
this apparent rarity of counterexamples surely means that something
interesting is going on, and we have found that this area suggests many
intriguing questions, a few, but not all, of which we discuss here.

We call a word \textit{symmetric} if it reads the same right to left as
left to right; e.g., \textit{ABBABBA} is symmetric, but \textit{ABABBA}
is not (in other contexts, the name ``palindromic'' is also used). To
simplify exposition, we shall often use exponents in the representation
of a word; e.g., the symmetric word above might have been written
$AB^2AB^2A$. We are principally concerned here with real symmetric
positive definite matrices, though in many cases the complex Hermitian
case is the same. We shall try to explicitly draw a distinction only
when it is important. We intend to exploit differences in the complex
Hermitian case in further work. Certain symmetries of a word do not
change the eigenvalues, and, since eigenvalues are our interest, we
shall freely use such symmetries and, often, only view two words as
distinct if they are not equivalent via the following transformations:

\begin{itemize}

\item 
\textit{Reversal}. Writing the letters of the words in reverse order.
This corresponds to transposition of the matrix product and thus does
not change eigenvalues.

\item 
\textit{Cyclic permutation}. Movement of the first letter of the word
to the end of the word. This can be realized as a similarity of the word
via the first letter and, thus, also does not change eigenvalues.

\item 
\textit{Interchange of $A$ and $B$}. This may change the eigenvalues of
a particular word, but, as $A$ and $B$ are both positive definite, it
does not change the possible eigenvalues.
\end{itemize}

Note that a symmetric word is one that is identical to its own reversal.
There are, for example, 20 words of length 6 with 3 $A$'s, but only 3
that are distinct up to the above symmetries: $ABABAB$, $A^3B^3$, and
$ABA^2B^2$.

Tangentially, we note that there is an algorithm for generating the
equivalence class, relative to the above symmetries, of a word of length
$L$ or determining the number of distinct equivalence classes among $N$
such words. Given a word $W$, another word $V$ lies in its equivalence
class if and only if $V$ is the result of $k$ cyclic permutations
$(0\leq k\leq L)$,
composed with (possibly) a reversal, composed with (possibly) an
interchange, applied to $W$. This gives an algorithm of order O$(NL)$.

Since a symmetric word may inductively be seen to be congruent
\cite[p.\ 223]{HJ} to either the center letter (if the length is odd) or to $I$ (if
the length is even), we have by Sylvester's law of inertia the
following.

\begin{lem}\label{lem1}
A symmetric word in two positive definite letters is positive definite
and, thus, has positive eigenvalues.
\end{lem}

It follows that any symmetric word gives an affirmative answer to
Question 3.

It has long been known \cite{HJ} that a product of two positive definite
matrices (e.g., the word $AB$) has positive eigenvalues and is
diagonalizable. We call a diagonalizable matrix with positive eigenvalues \textit{quasi-positive}
and record here a slightly more complete observation.

\begin{lem}\label{lem2}
The n-by-n matrix $Q$ is quasi-positive if and only if $Q=AB$, in which
$A$ and $B$ are positive definite. Moreover if $Q=SDS^{-1}$, with $D$ a
positive diagonal matrix, then all factorizations $AB$ of $Q$ into
positive definite matrices $A$ and $B$ are given by
$$
A=SES^*\quad\mbox{and}\quad B=S^{-1*}E^{-1}DS^{-1},
$$
in which $E$ is a positive definite matrix that commutes with $D$.
\end{lem}

\begin{proof}
If $Q=AB$, with $A$ and $B$ positive definite matrices, then $Q$ is
similar to $A^{-1/2}ABA^{1/2}=A^{1/2}BA^{1/2}$, which is congruent to
$B$ and, therefore, positive definite. Thus, $Q$ has positive
eigenvalues and is diagonalizable, as is so for a positive definite
matrix.

If $Q$ is quasi-positive, $Q=SDS^{-1}$, with $D$ positive diagonal, then
$Q=AB$, with $A=SES^*$ and $B=S^{-1*}E^{-1}DS^{-1}$ ($E$ is a positive
definite matrix commuting with $D$), both positive definite. Suppose
that $Q=AB$ is some other factorization into positive definite matrices. So
$B=A^{-1}Q$ is Hermitian. Then, $A^{-1}Q=Q^*A^{-1}$ or $AQ^*=QA$ or
$AS^{-1*}DS^*=SDS^{-1}A$, so that $S^{-1}AS^{-1*}D=DS^{-1}AS^{-1*}$.
Thus, $S^{-1}AS^{-1*}$ commutes with $D$; call $E=S^{-1}AS^{-1*}$, and
then $A=SES^*$. It follows that $E$ is Hermitian and positive definite,
as $A$ is. Now, $B=A^{-1}Q=S^{-1*}E^{-1}DS^{-1}$, which is positive
definite since $E^{-1}D$ is (because they commute).
\qquad\end{proof}

We now know that the nonsymmetric word $AB$ also positively answers
Question~3, but much more follows from Lemmas \ref{lem1} and \ref{lem2}.
We call a word \textit{nearly symmetric} if it is either symmetric or
the product (juxtaposition) of two symmetric words. It is an interesting
exercise that the nearly symmetric words are unchanged by the three
symmetries (i), (ii), and (iii). There is also a simple algorithm to
check for near symmetry: left to right, parse a given word after each
initial symmetric portion and check the remainder for symmetry (counting
the empty word as symmetric). We then have the following.

\begin{thm}\label{thm3}
Every nearly symmetric word in two positive definite letters has only
positive eigenvalues.
\end{thm}

\begin{proof}
The proof follows from Lemmas \ref{lem1} and \ref{lem2}.
\qquad\end{proof}

Are all words nearly symmetric? No, but all sufficiently short words
are.

\begin{thm}\label{thm4}
A word in which one of the letters appears at most twice is nearly
symmetric.
\end{thm}

\begin{proof}
Without loss of generality, we examine the situation in which $B$
appears at most twice. If a word contains only the letter $A$, the
result is trivial. If the letter $B$ appears only once, then the word
will be of the form $A^pBA^q(p,q\geq 0)$. If $p\geq q$, then we have
$A^pBA^q=A^{p-q}(A^qBA^q)$, and if $p\leq q$, we have
$A^pBA^q=(A^pBA^p)A^{q-p}$. In both cases, the word is nearly symmetric.
In the case of two $B$'s, the word can be written as $A^pBA^qBA^t
(p,q,t\geq 0)$, and so our word is one of the nearly symmetric words,
$(A^pBA^qBA^p)A^{t-p}$ or $A^{p-t}(A^tBA^qBA^t)$.
\qquad\end{proof}

In order to not be nearly symmetric then, a word must have length at
least 6 and 3 each of $A$ and $B$. Among the 3 such equivalence classes
of words of length 6, one is actually not nearly symmetric, $ABA^2B^2$,
and this shows that Theorem \ref{thm4} is best possible. This is the first
interesting word relative to Question 3, and we have the following
corollary.

\begin{cor}\label{cor5}
Every nearly symmetric word, and thus every word of length $<6$ has only
positive eigenvalues.
\end{cor}

An interesting question one can ask is how many nearly symmetric words
there are of a given length $L$. More importantly, what does the
fraction of nearly symmetric words to the total number of words approach
as $L$ goes to infinity? The result can be found in \cite{K}, and it
states that the number of nearly symmetric words of length $L$ is
O$(L\cdot 2^{(3/4)L})$. This gives us that the density of such words
approaches 0, and therefore, as $L$ goes to infinity, there is a pool of
potential negative answers to Questions 2 and 3 that ever increases in
relative frequency.

The situation is much simpler for 2-by-2 matrices, and we note (as does
Pierce \cite{P} and  Spitkovsky \cite{S}) the following.

\begin{fact}\label{fct1}
Both eigenvalues of any word in two $2$-by-$2$ positive definite
matrices are positive.
\end{fact}

\begin{proof}
We will actually show something stronger. Let $W$ be any finite product
of real positive powers of $A$ and $B$, in which $A$ and $B$ are 2-by-2
positive definite (complex) Hermitian matrices. (Here, we take principal
powers, so that $W$ is uniquely defined.) We first preprocess the word
as follows. Make one letter diagonal via uniform unitary similarity, and
then make the other letter entrywise nonnegative via a diagonal unitary
similarity. This does not change the first letter. Now, the word is
nonnegative (as it is clear from the spectral theorem that a positive
power of a nonnegative 2-by-2 positive definite matrix is nonnegative).
If it is diagonal, there is nothing more to do (the diagonal entries are
positive). If not, apply the Perron--Frobenius theorem (which says a
positive matrix must have a positive eigenvalue \cite[p.\ 503]{HJ}) and
the fact that the determinant is positive to show that the other
eigenvalue is positive as well.
\end{proof}

\begin{cor}\label{cor7}
The polynomial $p(t)$, defined by $p(t)=\mathrm{Tr}[(A+Bt)^m]$, has all
positive coefficients whenever $A$ and $B$ are $2$-by-$2$ positive definite
matrices.
\end{cor}

This all suggests that careful consideration of the word $ABA^2B^2$, or,
equivalently, $(BA)(BA)(AB)$, for 3-by-3 positive definite $A$ and $B$
is warranted. This is equivalent, by Lemma \ref{lem2}, to the study of
the expression $C^2C^{\mathrm{T}}$ for quasi-positive $C$. Since any
real matrix with real eigenvalues may be upper triangularized by
orthogonal similarity, it suffices to consider
$$
C=\left[\begin{array}{@{}ccc@{}}
a&x&z\\
0&b&y\\
0&0&c\end{array}\right]
$$
with $a,b,c>0$. If $a$, $b$, and $c$ are distinct, $C$ is diagonalizable
and thus quasi-positive. Using MAPLE, and with the assistance of Shaun
Fallat, it was found that $x,y,z$ and such $a,b,c$ may be found so that
Tr$(C^2C^{\mathrm{T}})<0$. Consistent with prior empirical experience,
choice of such $x,y,z$ and $a,b,c$ is delicate and falls in a very
narrow range. Resulting $A$ and $B$ (see Lemma \ref{lem2}) that exhibit
a negative answer to Question 2 (and, thus, 3) are, for example,
$$
A_1=\left[\begin{array}{@{}ccc@{}}
1&20&210\\
20&402&4240\\
210&4240&44903
\end{array}\right]\quad \mbox{and}\quad B_1=\left[\begin{array}{@{}ccc@{}}
36501&-3820&190\\
-3820&401&-20\\
190&-20&1
\end{array}\right].
$$
The extreme and reverse diagonal progressions are typical of such
examples. If the diagonal of one is ``flattened'' by orthogonal
similarity, the progression on the diagonal of the other becomes more
extreme.

We remark at this point that words giving a negative answer to Question 2
in the 3-by-3 case imply negative answers in the $n$-by-$n$ case for
$n>3$. This allows us to restrict our attention to the 3-by-3 positive
definite matrices. Simply direct sum a 3-by-3 example (giving a negative
trace) with a sufficiently small positive multiple of the identity to
get a larger example.

The idea of our first construction and some fortunate characteristics of
the constructed pair allow the identification of several infinite
classes of words giving negative answers to Questions 2 and 3. We
indicate some of these next.

1. Any positive integer power of a word that does not guarantee positive
eigenvalues also does not guarantee positive eigenvalues. For instance,
this shows that $BABAABBABAAB$ can have a nonpositive eigenvalue. This
is Theorem \ref{thm8} below.

2. Suppose a word can be written in terms of another word $T$ as
$T^k(T^*)^j$ for $k\neq j$. Furthermore, suppose $T=S_1S_2$ is a product
of two symmetric words $S_1$ and $S_2$. Then if the simultaneous word
equations
$$
\begin{array}{c}
S_1(A,B)=C,\\[3pt]
S_2(A,B)=D
\end{array}
$$
may be solved for positive definite $A$ and $B$ given positive definite
$C$ and $D$, then the original word can have negative trace. The first
nontrivial application of this technique is the first counterexample,
$(BA)^2AB$, in which $S_1=B$, $S_2=A$, $k=2$, and $j=1$. This result is
Theorem \ref{thm9} below.

3. Infinite classes involving single-letter length extension: this is a
nice application of sign analysis. Our first result is the following.

\begin{itemize}

\item[(a)]
\textit{The word, $ABA^2B^{2+k}$ with $k$ a nonnegative integer can have
negative trace.}

\end{itemize}

\begin{proof}
A direct computation with $A_1$ and $B_1$ from above gives us that
$$
(BABAAB)B=\left[\begin{array}{@{}ccc@{}}
-164679899&17226460&-856450\\
62354360&-6523192&324340\\
-5877450&614880&-30573
\end{array}\right]
$$
has sign pattern
$$
\left[\begin{array}{@{}ccc@{}}
-&+&-\\
+&-&+\\
-&+&-\end{array}
\right].
$$
Next, notice that $B_1$ has the sign pattern
$$
\left[\begin{array}{@{}ccc@{}}
+&-&+\\
-&+&-\\
+&-&+\end{array}\right]
$$
and that
$$
\left[\begin{array}{@{}ccc@{}}
-&+&-\\
+&-&+\\
-&+&-\end{array}\right]
\,
\left[\begin{array}{@{}ccc@{}}
+&-&+\\
-&+&-\\
+&-&+\end{array}\right]
$$
is
$$
\left[\begin{array}{@{}ccc@{}}
-&+&-\\
+&-&+\\
-&+&-\end{array}
\right]
$$
unambiguously.

Hence, multiplying the product $BABAABB$ by $B$ on the right any number
of times will preserve the negativity of the trace. Therefore,
$BABAABB\cdot B^k$ gives a negative answer to Question 2 for all
integers $k\geq 0$.
\end{proof}

Proofs using the same technique give us many infinite classes of
counterexamples, some of which we list below:

\begin{itemize}

\item[(b)]
$ABABAAB^k$, $k\geq 2$.

\item[(c)]
$ABBABAAB^k$, $k\geq 2$.

\item[(d)]
$ABAABBAAB^k$, $k\geq 2$.
\end{itemize}

4. Recall the two matrices $A_1$ and $B_1$ giving $(BA)(BA)(AB)$ a
negative trace. These matrices can also be used to prove that the words
$ABA^pB^q$, $ABBABA^pB^q$, and $ABABA^pB^q$ can have a negative trace
for all integers $p,q\geq 2$. Notice that (a), (b), and (c) above are
corollaries to this result. This is Theorem \ref{thm10} below.

We now present proofs of the three theorems mentioned above.

\begin{thm}\label{thm8}
Let $W$ be any word for which there are positive definite $A$ and $B$
such that $W(A,B)$ has an eigenvalue that is not positive. Then, for any
positive integer $k$, there are positive definite letters such that
$W^k$ has a nonpositive eigenvalue.
\end{thm}

\begin{proof}
Let $A$, $B$ be positive definite matrices that give $W$ a nonpositive
eigenvalue, and let $\lambda$ be such an eigenvalue. If
$k\in\{1,2,3,\ldots\}$, then an eigenvalue of $W(A,B)^k$ is $\lambda^k$.
If $\lambda^k$ is nonpositive, we are done, so the problem lies in the
possibility that $\lambda^k>0$. It will be necessary, therefore, in this
case to create a new pair of positive definite matrices $A'$ and $B'$
that give $W(A',B')^k$ a nonpositive eigenvalue.

We first offer a description of our approach before presenting the
details that follow. The idea is to parameterize a pair of positive
definite matrices in terms of a real variable $t$, $0\leq t \leq 1$, and
then examine the eigenvalues of the word $W^k$ evaluated at those
matrices. Using the continuity of eigenvalues on matrix entries, we then
show that $W(A(t),B(t))^k$ cannot have positive eigenvalues for all
$0\leq t\leq 1$.

Let $\lambda_A$ be the largest eigenvalue of $A$, and let $\lambda_B$ be
the largest eigenvalue of $B$. Define the following parameterization:
$$
A(t)=t\cdot(\lambda_AI-A)+A\quad\mbox{and}\quad B(t)=t\cdot(\lambda_BI-
B)+B\quad\mbox{for}\ 0\leq t\leq 1.
$$
We first note that $A(t)$ and $B(t)$ are positive definite for all such
$t$ since $(\lambda_AI-A)$ and $(\lambda_BI-B)$ are positive
semidefinite by a simple eigenanalysis. Next, notice that
$A(1)=\lambda_AI$ and $B(1)=\lambda_BI$, giving $W(A(1),B(1))$ positive
eigenvalues. Additionally, $A(0)=A$ and $B(0)=B$, which shows that
$W(A(0),B(0))$ has a nonpositive eigenvalue, by assumption. Since the
eigenvalues of a matrix depend continuously on its entries
\cite[p.\ 539]{HJ}, the eigenvalues of $W(A(t),B(t))$ also depend continuously on
$t$.

For $t\in [0,1]$, the spectrum of $W(A(t),B(t))$ cannot contain 0
because each product, $W(A(t),B(t))$, has positive determinant. Now, let
$$
\Gamma=\{t\in [0,1]\;|\;W(A(t),B(t))\ \mbox{has a positive spectrum}\}.
$$
Clearly, this set is not empty as $1\in \Gamma$, and it is not the
entire interval as $0\not\in \Gamma$. A straightforward continuity
argument also shows that $\Gamma$ is closed. Let $t_M$ be the greatest
lower bound of $\Gamma$, and notice that from above, $t_M\neq 0$ and
$t_M\in \Gamma$. As a result, the eigenvalues of $W(A(t_M),B(t_M))$ are
all positive. By continuity again, we can choose $t<t_M$ such that the
eigenvalues of $W(A(t),B(t))$ are as close to the eigenvalues of
$W(A(t_M),B(t_M))$ as we wish.

We are now ready to prove the theorem. Let $k$ be a positive integer.
By continuity, choose $t<t_M$
such that there is an eigenvalue, $\lambda$, of $W(A(t),B(t))$ with an
argument $\theta$ satisfying $-\pi/k<\theta<\pi/k$ (see Figure \ref{fig1}). This
guarantees that $\lambda^k$ cannot be real.
Our new pair $A(t),B(t)$ now proves the word $W^k$ can have nonpositive
eigenvalues.
\end{proof}


\begin{thm}\label{thm9}
If $j$ and $k$ are positive integers such that $j\neq k$, then there is
a real, quasi-positive matrix $T$ such that $T^k(T^*)^j$ has negative
trace.
\end{thm}

\begin{proof}
We first note that we can assume $k>j$, since if $k<j$, we examine
$[T^k(T^*)^j]^*$. We also assume without loss of generality that $T$ has 1
for an eigenvalue and it is the smallest eigenvalue of $T$.

Using Schur triangularization, we suppose
$$
T=\left[\begin{array}{@{}ccc@{}}
1&x&z\\
0&a&y\\
0&0&b\end{array}\right],
$$
with $x,y,z\in \Re$ and $b>a>1$.

Since it is necessary to compute powers of $T$, we note that
$$
T^k=\left[\begin{array}{@{}ccc@{}}
1&X_k&Z_k\\
0&a^k&Y_k\\
0&0&b^k\end{array}\right]=
\left[\begin{array}{@{}ccc@{}}
1&X_{k-1}&Z_{k-1}\\
0&a^{k-1}&Y_{k-1}\\
0&0&b^{k-1}\end{array}\right]\,
\left[\begin{array}{@{}ccc@{}}
1&x&z\\
0&a&y\\
0&0&b\end{array}\right],
$$
in which $X_k(Y_k$; $Z_k)$ is the 1,2 (2,3; 1,3) entry of $T^k$, $k>0$.

The above expression allows us to find formulae for the entries of $T^k$
by way of the following obvious recurrences:
$$
X_k=x+aX_{k-1};\quad Y_k=ya^{k-1}+bY_{k-1};\quad Z_k=z+yX_{k-1}+bZ_{k-
1}.
$$
An easy induction gives us that
$$
\begin{array}{c}
 X_k=x\frac{a^k-1}{a-1};\quad Y_k=y\frac{a^k-b^k}{a-b};\\[12pt]
 Z_k=xy\frac{1}{a-1}\cdot\left(\frac{a^k-b^k}{a-b}-b^{k-1}-
\frac{b^{k-1}-1}{b-1}\right)
+z\frac{b^k-1}{b-1}=xyC_k+zD_k,
\end{array}
$$
in which $C_k=\frac{1}{a-1}\cdot(\frac{a^k-b^k}{a-b}-b^{k-1}-
\frac{b^{k-1}-1}{b-1})$, $D_k=\frac{b^k-1}{b-1}$ depend only on $a$,
$b$, and $k$.

Thus, the trace of $T^k(T^*)^j$ can be computed explicitly in terms of
$x,y,z,a,b,k,j$. It is
$$
\begin{array}{l}

\mathrm{Tr}[T^k(T^*)^j]=\mathrm{Tr}\left[\left(\begin{array}{@{}ccc@{}}
1&X_k&Z_k\\
0&a^k&Y_k\\
0&0&b^k\end{array}\right)\left(\begin{array}{@{}ccc@{}}
1&0&0\\
X_j&a^j&0\\
Z_j&Y_j&b^j\end{array}\right)\right]\\[24pt]
\qquad = (1+X_kX_j+Z_kZ_j)+(a^{k+j}+Y_kY_j)+b^{k+j}\\[12pt]
\qquad\displaystyle = 1+a^{k+j}+b^{k+j}+x^2\frac{a^k-1}{a-1}\cdot\frac{a^j-1}{a-1}+
y^2\frac{a^k-b^k}{a-b}\cdot \frac{a^j-b^j}{a-b}\\[12pt]
\qquad\qquad +x^2y^2C_kC_j+xyz(C_kD_j+C_jD_k)+z^2D_kD_j.
\end{array}
$$
Fix $a,b>1$ and set $y=x$. Now, view Tr$[T^k(T^*)^j]$ as a quadratic
polynomial in $z$. For this polynomial to take on negative values, it is
necessary and sufficient for its discriminant to be positive. This
discriminant is a quartic polynomial in $x$; therefore, if we can show
that its leading coefficient is always positive, this will demonstrate
that for large enough values of $x$, the discriminant will also be
positive. The coefficient of $x^4$ in this discriminant is
$$
\begin{array}{c}
(C_kD_j+C_jD_k)^2-4D_kD_j(C_kC_j)\\[12pt]
=C^2_kD^2_j+2C_kC_jD_kD_j+C^2_jD^2_k-4C_kC_jD_kD_j=(C_kD_j-C_jD_k)^2.
\end{array}
$$
When $k=j$, the expression above is 0, so it is necessary to prove that
whenever $k\neq j$, $C_kD_j\neq C_jD_k$. Examining $C_kD_j-C_jD_k$, this
is equivalent to proving that
$$
a^jb^{k+1}-a^k+b^k+a^kb-b^{k+1}+a^kb^j-a^kb^{j+1}-a^jb^k-b^j+a^j-
a^jb+b^{j+1}
$$
is never zero unless $k=j$. Factoring out $(b-1)$, we need only prove
that
$$
f(a,b)=a^jb^k+a^k-b^k-a^kb^j+b^j-a^j
$$
is never zero unless $k=j$. Examine the following polynomial in $x$:
$$
g(x)=f(a,x)=x^k(a^j-1)+x^j(1-a^k)+a^k-a^j.
$$
It is easy to see that $g(1)=0$ and $g(a)=0$. From Descartes's rule of
signs, it is clear (since $a>1$) that $g$ has either 0 or 2 positive
real roots. Since $a$ and 1 are two such roots, $g$ has no more positive
ones. Hence, $g(b)\neq 0$ for $b\neq 1$, $a$.

This concludes the proof that $T^k(T^*)^j$ will have negative trace for
some quasi-positive matrix $T$. Note that a description of all 3-by-3
quasi-positive $T$ that give $T^k(T^*)^j$ a negative trace is implicit
in the proof.
\end{proof}

Our first corollary to this theorem is that the word $(BA)^2AB$ gives a
negative answer to Question 2; but moreover, we also now have a
description of all 3-by-3 positive definite $A$ and $B$ that give
$(BA)^2AB$ a negative trace. Theorem \ref{thm9} describes all 3-by-3
quasi-positive matrices $T$ that give $T^2T^*$ a negative trace, and
hence all positive definite matrices $A$ and $B$ are given by $T=BA$
from Lemma \ref{lem2}.

We now prove the following.

\begin{thm}\label{thm10}
For integers $p,q\geq 2$ and the word $W=ABA^pB^q$, there exist positive
definite matrices $A$ and $B$ such that $W(A,B)$ has a negative trace.
\end{thm}

\begin{proof}
We first record a few preliminaries.

Let $F(p,q)=\mathrm{Tr}[ABA^pB^q]=\mathrm{Tr}[BAB^qA^p]$ be the desired
trace of the word $W$. Now, suppose $A=U^*DU$ and $B=V^*EV$ are fixed
positive definite matrices with $U,V$ (real) orthogonal, and let
$D=\text{diag}(a,b,c)$, $E=\text{diag}(r,s,t)$, $a,b,c,r,s,t>0$. Then we can write
$$
F(p,q)=\mathrm{Tr}[UBAB^qU^*D^p]=\mathrm{Tr}[VABA^pV^*E^q].
$$
From these two expressions, it is clear that
\begin{eqnarray}
F(p,q)&=&g_1(q)a^p+g_2(q)b^p+g_3(q)c^p,\label{eq1}\\[3pt]
F(p,q)&=&h_1(p)r^q+h_2(p)s^q+h_3(p)t^q,\label{eq2}
\end{eqnarray}
where $g_i(q)$, $h_i(p)$ are linear functions in $r^q,s^q,t^q$ and
$a^p,b^p,c^p$, respectively. Equations (\ref{eq1}) and (\ref{eq2}) can
be viewed as a generalization of the well-known expression for computing
Fibonacci numbers. In fact, these equations imply the recurrence
relations
\begin{eqnarray} 
\quad F(p,q)&=&(a+b+c)F(p-1,q)-(ab+bc+ac)F(p-2,q)+(abc)F(p-
3,q),\label{eq3}\\[3pt]
\quad F(p,q)&=&(r+s+t)F(p,q-1)-(rs+rt+st)F(p,q-2)+(rst)F(p,
q-3).\label{eq4}
\end{eqnarray}
We are now ready to prove the result. It turns out that $A_1$ and $B_1$
(as described above) will prove the claim
$$
A_1=\left[\begin{array}{@{}ccc@{}}
1&20&210\\
20&402&4240\\
210&4240&44903\end{array}\right],
\qquad B_1=\left[\begin{array}{@{}ccc@{}}
36501&-3820&190\\
-3820&401&-20\\
190&-20&1
\end{array}\right].
$$
The values of $(a+b+c)$, $(ab+bc+ac)$, $(abc)$, $(r+s+t)$, $(rs+st+rt)$,
and $(rst)$ are obtained from the characteristic polynomials of $A$ and
$B$. These polynomials are easy to compute as $P_A(t)=t^3-
45306t^2+74211t-6$ and $P_B(t)=t^3-36903t^2+44903t-1$. Therefore,
(\ref{eq3}) and (\ref{eq4}) become
\begin{eqnarray} 
F(p,q)&=&45306\cdot F(p-1,q)-74211\cdot F(p-2,q)+6\cdot F(p-3,q),\label{eq5}\\[3pt]
F(p,q)&=&36903\cdot F(p,q-1)-44903\cdot F(p,q-2)+F(p,q-3).\label{eq6}
\end{eqnarray}
To prove the theorem, we must show that $F(p,q)<0$ for all $p,q\geq 2$.
First notice that for the base cases of $2\leq p$, $q\leq 4$, we have
that $F(p,q)$ are given by the following table:
\vspace*{12pt}

\begin{center}
\noindent{\footnotesize%
\begin{tabular}{|l|c|c|c|}
\hline
&\boldmath$q=2$\rule{0pt}{9pt}&\boldmath$q=3$&\boldmath$q=4$\\[1pt]
\hline
$p=2$\rule{0pt}{9pt}&$-$3164&$-$171233664&$-$6318893781764\\[1pt]
\hline
$p=3$\rule{0pt}{9pt}&$-$219049002&$-$10537988104302&$-$388873536893369802\\[1pt]
\hline
$p=4$\rule{0pt}{9pt}&$-$9923997300324&$-$477421308542380824
&$-$17617832833924812095724\\[2pt]
\hline
\end{tabular}}
\end{center}
\vspace*{12pt}

\noindent To prove the result using the recurrences above, we will
invoke induction and prove something stronger. Namely, we claim that for
all $p,q\geq 2$, $F(p,q)<0$ and also the following inequalities hold:
$$
\begin{array}{c}
F(p,q)<10\cdot F(p-1,q)\quad \mbox{for}\ p>2,\\[3pt]
F(p,q)<10\cdot F(p,q-1)\quad \mbox{for}\ q>2.
\end{array}
$$
Suppose the result is true for all $2\leq p,q<N$ (from the table above,
we can also suppose $N\geq 5$); then we want to show it true for $2\leq p$,
$q\leq N$. For $2\leq p,q<N$, examine $F(N,q), F(p,N)$, and $F(N,N)$.
From (\ref{eq5}) and (\ref{eq6}), we have
\begin{equation}\label{eq7}
\begin{array}{r@{\;}l}
\quad\,  F(N,q)& = 45306\cdot F(N-1,q)-74211\cdot F(N-2,q)+6\cdot F(N-3,q)\\[3pt]
& < 45306\cdot F(N-1,q)-7421.1\cdot F(N-1),q)\\[3pt]
& = 37884.9\cdot F(N-1,q)<10\cdot F(N-1,q),
\end{array}
\end{equation}
\vspace*{-12pt}
\begin{equation}\label{eq8}
\begin{array}{r@{\;}l}
F(p,N)& = 36903\cdot F(p,N-1)-44903\cdot F(p,N-2)+F(p,N-3)\\[3pt]
& < 36903\cdot F(p,N-1)-4490.3\cdot F(p,N-1)\\[3pt]
& = 32412.7\cdot F(p,N-1)<10\cdot F(p,N-1).
\end{array}
\end{equation}
But to complete the induction, we must also show that $F(N,N)<10\cdot
F(N-1,N)$ and $F(N,N)<10\cdot F(N,N-1)$. Substituting (\ref{eq6}) into
the right-hand side of (\ref{eq5}) with $p=N$, $q=N$, we have
\begin{equation}\label{eq9}
\begin{array}{@{\extracolsep{-.5pc}}rcl}
F(N,N)& = & 1671927318\cdot F(N-1,N-1)-2034375318\cdot F(N-1,N-2)\\[3pt]
&&+\ 45306\cdot F(N-1,N-3)\\[3pt]
&&  -\ 2738608533\cdot F(N-2,N-1)+3332296533\cdot F(N-2,N-2)\\[3pt]
&& -\ 74211\cdot F(N-2,N-3)\\[3pt]
&& +\ 221418\cdot F(N-3, N-1)-269418\cdot F(N-3, N-2)\\[3pt]
&& +\ 6\cdot F(N-3, N-3)\\[3pt]
& < &1671927318\cdot F(N-1,N-1)-203437531.8\cdot F(N-1,N-1)\\[3pt]
&&  -\ 273860853.3\cdot F(N-1,N-1)-74.211\cdot F(N-1,N-1)\\[3pt]
&& -\ 269.418\cdot F(N-1,N-1)\\[3pt]
& =&1194628589.271\cdot F(N-1,N-1).
\end{array}
\end{equation}
But from (\ref{eq8}) with $p=N-1$, we have
$$
\begin{array}{r@{\;}l}
36903\cdot F(N-1,N-1)& = F(N-1,N)+44903\cdot F(N-1,N-2)-F(N-1,N-
3)\\[3pt]
& < F(N-1,N)-(1/100)\cdot F(N-1,N-1).
\end{array}
$$
Therefore, $36903.01\cdot F(N-1,N-1)<F(N-1,N)$, which gives us easily
(from (\ref{eq9})) that
$$
F(N,N)<1194628589.271\cdot F(N-1,N-1)<10\cdot F(N-1,N).
$$
To arrive at $F(N,N)<10\cdot F(N,N-1)$, we perform the same
examination, this time with (\ref{eq7}):
$$
F(N,N-1)=45306\cdot F(N-1,N-1)-74211\cdot F(N-2,N-1)+6\cdot F(N-3,N-1),
$$
giving us the inequality $45306.06\cdot F(N-1,N-1)<F(N,N-1)$.

So again, from (\ref{eq9}), we see that $F(N,N)<10\cdot F(N,N-1)$. This
completes the induction and shows that for all $p,q\geq 2$, $F(p,q)<0$.
The proof also bounds the growth from below, but the factor of 10 is
obviously not the best possible.
\end{proof}

At this point, we should remark that the proof for Theorem \ref{thm10}
above could be generalized to a certain extent. Namely, suppose $W$ is a
word that can be written as $W_1A^pW_2B^q$ for some words $W_1,W_2$ in
$A$ and $B$. Then, $A_1$ and $B_1$ give this word negative trace for all
integers $p,q\geq 2$ provided that for the base cases of $2\leq p$,
$q\leq 4$,
$$
F(p,q)<0;\quad F(p,q)<10\cdot F(p-1,q);\quad\mbox{and}\ F(p,q)<10\cdot
F(p,q-1).
$$
As an example, a calculation gives us that for the word $W=ABABA^pB^q$
the first 9 values of $F(p,q)$ are given by\footnote{While values are
integers, they are shown only to the first 15 significant digits.}
\vspace*{12pt}

\begin{center}
\noindent{\footnotesize
\begin{tabular}{|l|c|c|c|}
\hline
&\boldmath$q=2$\rule{0pt}{9pt}&\boldmath$q=3$&\boldmath$q=4$\\[1pt]
\hline
$p=2$\rule{0pt}{9pt}&$-$32302&$-$1319655482&$-$48697748014592\\[1pt]
\hline
$p=3$\rule{0pt}{9pt}&$-$1748875224&$-$70292975950848&$-
$2.59394099689082e$+$018\\[1pt]
\hline
$p=4$\rule{0pt}{9pt}&$-$79232137801728&$-$3.18459541653658e$+$018&$-
$1.17517468821039e$+$023\\[1pt]
\hline
\end{tabular}}
\end{center}
\vspace*{12pt}

\noindent The word $W=ABBABA^pB^q$ also satisfies the base case
conditions as the $F(p,q)$ are
\vspace*{12pt}

\begin{center}
\noindent{\footnotesize
\begin{tabular}{|l|c|c|c|}
\hline
&\boldmath$q=2$\rule{0pt}{9pt}&\boldmath$q=3$&\boldmath$q=4$\\[1pt]
\hline
$p=2$\rule{0pt}{9pt}&$-$222790424&$-$10720038844524&$-
$3.95591587257758e$+$017\\[1pt]
\hline
$p=3$\rule{0pt}{9pt}&$-$10103386100406&$-$4.86025787321779e$+$017&$-
$1.79353558546523e$+$022\\[1pt]
\hline
$p=4$\rule{0pt}{9pt}&$-$4.57727477164142e$+$017&$-
$2.20190887755731e$+$022&$-$8.12549875102683e$+$026\\[1pt]
\hline
\end{tabular}}
\end{center}
\vspace*{12pt}

\noindent It should now be clear that we conjecture the following.

\begin{conj}\label{con4}
{\rm A word has positive trace for every pair of positive
definite letters if and only if the word is nearly symmetric.}
\end{conj}

Using the results and ideas we have discussed, it is possible to verify
this conjecture for words of lengths less than 11. Before listing these
results, we remark on how to find specific $A$ and $B$ for which a word
has negative trace. One difficulty is how to view the set of positive
definite matrices $A$ and $B$. We explain a helpful parametric approach
for the sample word $BAABBAAA$ and the generalization will be clear.
Notice that we do not yet know that this word can have a negative trace
using any of the methods thus far.

First set $Q=AB$, and recall that all solutions $A$, $B$ to such an
equation are given by Lemma \ref{lem2} as $Q=SDS^{-1}$, $A=SES^*$,
$B=S^{-1*}E^{-1}DS^{-1}$, in which $D$ is a positive diagonal matrix,
and $E$ is a positive definite matrix commuting with $D$. For
simplicity, we seek a positive diagonal $E$. Using these substitutions
and some simplification, our original word has the same eigenvalues as
the following expression: $DPDP^{-1}DPEPE$, in which $P=S^*S$.

\begin{table}[t!]
\footnotesize
\caption{All words that are not nearly symmetric of length $< 11$ admit
negative trace.}
\label{tab1}
\begin{center}
\begin{tabular}{|l|l|}
\hline
$AABABB$\rule{0pt}{9pt}&Original solution $A_1, B_1$ using $C^2C^T$\\[1pt]
\hline
$AAABABB$\rule{0pt}{9pt}&Theorem \ref{thm10}\\[1pt]
\hline
$AAAABABB$\rule{0pt}{9pt}&Theorem \ref{thm10}\\[1pt]
\hline
$AAABAABB$\rule{0pt}{9pt}&Using $A_2,B_2$\\[1pt]
\hline
$AAABABBB$\rule{0pt}{9pt}&Theorem \ref{thm10}\\[1pt]
\hline
$AABABABB$\rule{0pt}{9pt}&Equivalent to $(AB)^3BA$\\[1pt]
\hline
$AAAAABABB$\rule{0pt}{9pt}&Theorem \ref{thm10}\\[1pt]
\hline
$AAAABAABB$\rule{0pt}{9pt}&Equivalent to $(A^2B)^2BA^2$\\[1pt]
\hline
$AAAABABBB$\rule{0pt}{9pt}&Theorem \ref{thm10}\\[1pt]
\hline
$AAABAABAB$\rule{0pt}{9pt}&Using $A_3$, $B_3$ produced by the technique above\\[1pt]
\hline
$AAABAABBB$\rule{0pt}{9pt}&Using $A_2$, $B_2$\\[1pt]
\hline
$AAABABABB$\rule{0pt}{9pt}&Theorem \ref{thm10}\\[1pt]
\hline
$AABAABABB$\rule{0pt}{9pt}&Theorem \ref{thm10}\\[1pt]
\hline
$AAAAAABABB$\rule{0pt}{9pt}&Theorem \ref{thm10}\\[1pt]
\hline
$AAAAABAABB$\rule{0pt}{9pt}&Using $A_2$, $B_2$\\[1pt]
\hline
$AAAAABABBB$\rule{0pt}{9pt}&Theorem \ref{thm10}\\[1pt]
\hline
$AAAABAAABB$\rule{0pt}{9pt}&Using $A_4$, $B_4$ produced by the technique above\\[1pt]
\hline
$AAAABAABAB$\rule{0pt}{9pt}&Using $A_3$, $B_3$\\[1pt]
\hline
$AAAABAABBB$\rule{0pt}{9pt}&Using $A_2$, $B_2$\\[1pt]
\hline
$AAAABABABB$\rule{0pt}{9pt}&Theorem \ref{thm10}\\[1pt]
\hline
$AAAABABBBB$\rule{0pt}{9pt}&Theorem \ref{thm10}\\[1pt]
\hline
$AAAABBABBB$\rule{0pt}{9pt}&Using $A_2$, $B_2$ (interchanging $A$ and $B$)\\[1pt]
\hline
$AAABAABABB$\rule{0pt}{9pt}&Theorem \ref{thm10}\\[1pt]
\hline
$AAABAABBAB$\rule{0pt}{9pt}&Using $A_1$, $B_1$\\[1pt]
\hline
$AAABABAABB$\rule{0pt}{9pt}&Using $A_2$, $B_2$\\[1pt]
\hline
$AAABABABBB$\rule{0pt}{9pt}&Using $A_2$, $B_2$\\[1pt]
\hline
$AAABABBABB$\rule{0pt}{9pt}&Theorem \ref{thm10}\\[1pt]
\hline
$AAABBAABBB$\rule{0pt}{9pt}&Using $A_5$, $B_5$ produced by the technique above\\[1pt]
\hline
$AABABABABB$\rule{0pt}{9pt}&Equivalent to $(AB)^4BA$\\[1pt]
\hline
$AABABABBAB$\rule{0pt}{9pt}&Equivalent to $(AB)^3(BA)^2$\\[1pt]
\hline
$AABABBAABB$\rule{0pt}{9pt}&(d)\\[1pt]
\hline
\end{tabular}
\end{center}
\end{table}

Next, fix a positive definite matrix $P$ and view the positive diagonal
matrices $D$ and $E$ parametrically, hoping now to minimize the trace of
the product above. These minimizations are easier to perform because now
we have a simple parametric description of positive definite pairs.
Notice that it is not necessary to find $A$ and $B$ to show that they
exist and give the word a negative trace. However, it is useful to have
explicit examples, as they may be later used to show that other (not
nearly symmetric) words admit negative trace. After finding $D$, $E$,
and $P$, we recover these letters from the equations $S^*S=P$,
$A=SES^*$, $B=S^{-1*}E^{-1}DS^{-1}$. An example solution found using this
technique for the word $BAABBAAA$ is given by
$$
\begin{array}{r@{\;}l}
A_2& = \left[\begin{array}{@{}ccc@{}}
4351/479&4856/399&18421/62\\
4856/399&16073/64&3784/21\\
18421/62&3784/21&89917/9\end{array}\right],\\[18pt]
B_2&  = \left[\begin{array}{@{}ccc@{}}
2461/149&-297/641&-757/1569\\
-297/641&179/6146&50/3767\\
-757/1569&50/3767&269/19081\end{array}\right].
\end{array}
$$
It is easily verified that the trace of the word $BAABBAAA$ is a
negative rational number given approximately by
Tr$(BAABBAAA)\approx -143370.8471$.

In Table \ref{tab1} we list all the equivalence classes of words that
are not nearly symmetric and are of length less than 11. Next to each
word, we describe the method of finding the $A$ and $B$ that proves they
can have a negative trace.

\section{Acknowledgment}
The first author would like to acknowledge the pleasant and sometimes
useful conversations with several mathematicians about this problem---in
particular David Yopp and Tom Laffey.


\end{document}